\documentclass[12pt,twoside]{article}

\usepackage{amsmath,amsthm,amsfonts,amssymb}
\usepackage{hyperref}

\def\proof{{\bf Proof. }}

\date{} 

\begin{document} 

\centerline{} 

\centerline{} 

%%%%%%%%%%%%%%%%%%%%%%%%%%%%%%%%%%%%%%%%%%%%%%%%%%%%%%%%%%%%%%%%%%%%%%%%%%%%%%%%%%%%%%%%%%%%%%%%%%%%%%%%%%%%%%%%%%%%%%%%%%%%%%%%%%%%%%%%%%%%%%%%%%%%%%%%
\centerline {\Large{\bf On the Primes in the Interval $[3n, 4n]$}} 
%%%%%%%%%%%%%%%%%%%%%%%%%%%%%%%%%%%%%%%%%%%%%%%%%%%%%%%%%%%%%%%%%%%%%%%%%%%%%%%%%%%%%%%%%%%%%%%%%%%%%%%%%%%%%%%%%%%%%%%%%%%%%%%%%%%%%%%%%%%%%%%%%%%%%%%%

\centerline{} 

\centerline{\bf {Andy Loo}} 

\centerline{St. Paul's Co-educational College, Hong Kong} 

\centerline{} 

\newtheorem{Theorem}{\quad Theorem}[section] 

\newtheorem{Definition}[Theorem]{\quad Definition} 

\newtheorem{Corollary}[Theorem]{\quad Corollary} 

\newtheorem{Lemma}[Theorem]{\quad Lemma} 

\newtheorem{Example}[Theorem]{\quad Example} 

%%%%%%%%%%%%%%%%%%%%%%%%%%%%%%%%%%%%%%%%%%%%%%%%%%%%%%%%%%%%%%%%%%%%%%%%%%%%%%%%%%%%%%%%%%%%%%%%%%%%%%%%%%%%%%%%%%%%%%%%%%%%%%%%%%%%%%%%%%%%%%%%%%%%%%%%
%%%%%%%%%%%%%%%%%%%%%%%%%%%%%%%%%%%%%%%%%%%%%%%%%%%%%%%%%%%%%%%%%%%%%%%%%%%%%%%%%%%%%%%%%%%%%%%%%%%%%%%%%%%%%%%%%%%%%%%%%%%%%%%%%%%%%%%%%%%%%%%%%%%%%%%%

\begin{abstract} 
For the old question whether there is always a prime in the interval $[kn, (k + 1)n]$ or not,
the famous Bertrand's postulate gave an affirmative answer for $k = 1$. It was first proved by P.L. Chebyshev in 1850,
and an elegant elementary proof was given by P. Erd\H{o}s in 1932 (reproduced in \cite[pp.\,171-173]{erdos}).
M. El Bachraoui used elementary techniques to prove the case $k = 2$ in 2006 \cite{bachraoui}.
This paper gives a proof of the case $k = 3$, again without using the prime number theorem or any deep analytic result.
In addition we give a lower bound for the number of primes in the interval $[3n, 4n]$,
which shows that as $n$ tends to infinity, the number of primes in the interval $[3n, 4n]$ goes to infinity.
\end{abstract} 

%%%%%%%%%%%%%%%%%%%%%%%%%%%%%%%%%%%%%%%%%%%%%%%%%%%%%%%%%%%%%%%%%%%%%%%%%%%%%%%%%%%%%%%%%%%%%%%%%%%%%%%%%%%%%%%%%%%%%%%%%%%%%%%%%%%%%%%%%%%%%%%%%%%%%%%%
%%%%%%%%%%%%%%%%%%%%%%%%%%%%%%%%%%%%%%%%%%%%%%%%%%%%%%%%%%%%%%%%%%%%%%%%%%%%%%%%%%%%%%%%%%%%%%%%%%%%%%%%%%%%%%%%%%%%%%%%%%%%%%%%%%%%%%%%%%%%%%%%%%%%%%%%



\setcounter{section}{-1} 

%%%%%%%%%%%%%%%%%%%%%%%%%%%%%%%%%%%%%%%%%%%%%%%%%%%%%%%%%%%%%%%%%%%%%%%%%%%%%%%%%%%%%%%%%%%%%%%%%%%%%%%%%%%%%%%%%%%%%%%%%%%%%%%%%%%%%%%%%%%%%%%%%%%%%%%%
%%%%%%%%%%%%%%%%%%%%%%%%%%%%%%%%%%%%%%%%%%%%%%%%%%%%%%%%%%%%%%%%%%%%%%%%%%%%%%%%%%%%%%%%%%%%%%%%%%%%%%%%%%%%%%%%%%%%%%%%%%%%%%%%%%%%%%%%%%%%%%%%%%%%%%%%

\section{Notations} 

%%%%%%%%%%%%%%%%%%%%%%%%%%%%%%%%%%%%%%%%%%%%%%%%%%%%%%%%%%%%%%%%%%%%%%%%%%%%%%%%%%%%%%%%%%%%%%%%%%%%%%%%%%%%%%%%%%%%%%%%%%%%%%%%%%%%%%%%%%%%%%%%%%%%%%%%
%%%%%%%%%%%%%%%%%%%%%%%%%%%%%%%%%%%%%%%%%%%%%%%%%%%%%%%%%%%%%%%%%%%%%%%%%%%%%%%%%%%%%%%%%%%%%%%%%%%%%%%%%%%%%%%%%%%%%%%%%%%%%%%%%%%%%%%%%%%%%%%%%%%%%%%%

Throughout this paper, we let $n$ run through the positive integers and $p$ run through the primes.
We also let $\pi(n)$ be the prime counting function, which counts the number of primes not exceeding $n$.
Further define
	$$f(x) = \sqrt{2\pi}x^{x + \frac{1}{2}}e^{-x}e^{\frac{1}{12x}}$$
and
	$$g(x) = \sqrt{2\pi}x^{x + \frac{1}{2}}e^{-x}e^{\frac{1}{12x + 1}}.$$ 

%%%%%%%%%%%%%%%%%%%%%%%%%%%%%%%%%%%%%%%%%%%%%%%%%%%%%%%%%%%%%%%%%%%%%%%%%%%%%%%%%%%%%%%%%%%%%%%%%%%%%%%%%%%%%%%%%%%%%%%%%%%%%%%%%%%%%%%%%%%%%%%%%%%%%%%%
%%%%%%%%%%%%%%%%%%%%%%%%%%%%%%%%%%%%%%%%%%%%%%%%%%%%%%%%%%%%%%%%%%%%%%%%%%%%%%%%%%%%%%%%%%%%%%%%%%%%%%%%%%%%%%%%%%%%%%%%%%%%%%%%%%%%%%%%%%%%%%%%%%%%%%%%

\section{Lemmas} 

%%%%%%%%%%%%%%%%%%%%%%%%%%%%%%%%%%%%%%%%%%%%%%%%%%%%%%%%%%%%%%%%%%%%%%%%%%%%%%%%%%%%%%%%%%%%%%%%%%%%%%%%%%%%%%%%%%%%%%%%%%%%%%%%%%%%%%%%%%%%%%%%%%%%%%%%
%%%%%%%%%%%%%%%%%%%%%%%%%%%%%%%%%%%%%%%%%%%%%%%%%%%%%%%%%%%%%%%%%%%%%%%%%%%%%%%%%%%%%%%%%%%%%%%%%%%%%%%%%%%%%%%%%%%%%%%%%%%%%%%%%%%%%%%%%%%%%%%%%%%%%%%%

\begin{Lemma}\label{Lemma:oneone}
If $n \ge 8$, then
	$$\pi(n) \le \frac{n}{2}.$$
\end{Lemma}

\proof This is trivial since $1, 9$ and all even positive integers are not prime.

%%%%%%%%%%%%%%%%%%%%%%%%%%%%%%%%%%%%%%%%%%%%%%%%%%%%%%%%%%%%%%%%%%%%%%%%%%%%%%%%%%%%%%%%%%%%%%%%%%%%%%%%%%%%%%%%%%%%%%%%%%%%%%%%%%%%%%%%%%%%%%%%%%%%%%%%

\begin{Lemma}\label{Lemma:onetwo}
If $x$ is a positive real number, then
	$$\prod_{p \le x} p \le 4^x.$$
\end{Lemma}

\proof See \cite[pp.\,167-168]{erdos}.

%%%%%%%%%%%%%%%%%%%%%%%%%%%%%%%%%%%%%%%%%%%%%%%%%%%%%%%%%%%%%%%%%%%%%%%%%%%%%%%%%%%%%%%%%%%%%%%%%%%%%%%%%%%%%%%%%%%%%%%%%%%%%%%%%%%%%%%%%%%%%%%%%%%%%%%%

\begin{Lemma}\label{Lemma:onethree}
We have
	$$g(n) < n! < f(n)$$
\end{Lemma}

\proof See \cite{robbins}.

%%%%%%%%%%%%%%%%%%%%%%%%%%%%%%%%%%%%%%%%%%%%%%%%%%%%%%%%%%%%%%%%%%%%%%%%%%%%%%%%%%%%%%%%%%%%%%%%%%%%%%%%%%%%%%%%%%%%%%%%%%%%%%%%%%%%%%%%%%%%%%%%%%%%%%%%

\begin{Lemma}\label{Lemma:onefour}
For a fixed constant $c \ge \dfrac{1}{12}$, define the function
	$$h_1(x) = \frac{f(x + c)}{g(c)g(x)}.$$
Then for $x \ge \dfrac{1}{2}$, $h_1(x)$ is increasing. 
\end{Lemma}

\proof It suffices to prove that the function
	$$H_1(x) = (x + c)^{x + c + \frac{1}{2}}x^{-x - \frac{1}{2}}e^{\frac{1}{12(x + c)}- \frac{1}{12x + 1}}$$
is increasing for $x > \frac{1}{2}$. Indeed, we have
	$$H_1^\prime(x) = H_1(x)\Biggl(\biggl(\frac{1}{2(x + c)} - \frac{1}{12(x + c)^2} + \ln(x + c)\biggr)
	                  - \biggl( \frac{1}{2x} - \frac{1}{12\biggl(x + \dfrac{1}{12}\biggr)^2} + \ln x \biggr)\Biggr),$$
where $H_1(x) > 0$. Let
	$$F_1(x) = \frac{1}{2x} + \ln x.$$
As $-\dfrac{1}{12(x + c)^2} \ge \dfrac{1}{12\biggl(x + \dfrac{1}{12}\biggr)^2}$, it suffices to prove that $F_1(x)$ is increasing, so that
	$$\biggl(\frac{1}{2(x + c)} + \ln(x + c)\biggr) - \biggl(\frac{1}{2x} + \ln x\biggr) \ge 0.$$
We actually have
	$$F^\prime_1(x) = -\frac{1}{2x^2} + \frac{1}{x} = \frac{2x - 1}{2x^2},$$
which must be non-negative for all $x \ge \dfrac{1}{2}$. Therefore, the desired result follows.

\begin{Lemma}\label{Lemma:onefive}
For a fixed positive constant $c$, define the function
	$$h_2(x) = \frac{f(c)}{g(x)g(c - x)}.$$
Then when $\dfrac{1}{2} \le x < \dfrac{c}{2}$, $h^\prime_2(x) > 0$; when $x = \dfrac{c}{2}$, $h^\prime_2(x) = 0$ 
and when $\dfrac{c}{2} < x \le c - \dfrac{1}{2}$, $h^\prime_2(x) < 0$.
\end{Lemma}

\proof It suffices to prove that the function
	$$H_2(x) = x^{x + \frac{1}{2}}(c - x)^{c - x + \frac{1}{2}}e^{\frac{1}{12x + 1} + \frac{1}{12(c - x) + 1}}$$
has the following property: when $\dfrac{1}{2} \le x < \dfrac{c}{2}$, $H^\prime_2(x) < 0$; when $x = \dfrac{c}{2}$, $H^\prime_2(x) = 0$ 
and when $\dfrac{c}{2} < x \le c - \dfrac{1}{2}$, $H^\prime_2(x) > 0$. Indeed, we have
	$$H^\prime_2(x) = x^{x + \frac{1}{2}}(c - x)^{c - x + \frac{1}{2}}e^{\frac{1}{12x + 1} - \frac{1}{12(c - x) + 1}}(F_2(x) - F_2(c - x)),$$
where
	$$F_2(x) = \frac{12}{(12(c - x) + 1)^2} + \frac{1}{2x} + \ln x.$$
Clearly, for $\dfrac{1}{2} \le x < c$,
	$$x^{x + \frac{1}{2}}(c - x)^{c - x + \frac{1}{2}}e^{\frac{1}{12x + 1} - \frac{1}{12(c - x) + 1}} > 0.$$
Next, we actually have
	$$F^\prime_2(x) = \frac{288}{(12(c - x) + 1)^3} + \frac{2x - 1}{2x^2},$$
which must be positive for all $\dfrac{1}{2} \le x < c$. 
Thus whenever $\dfrac{1}{2} \le x \le c- \dfrac{1}{2}$, $F_2(x)$ is increasing while $F_2(c - x)$ is decreasing, implying that there is at most
one value of $x$ with $\dfrac{1}{2} \le x \le c- \dfrac{1}{2}$ satisfying $F_2(x) = F_2(c - x)$.
It is clear that $x = \dfrac{c}{2}$ is such a value. It follows that when $\dfrac{1}{2} \le x < \dfrac{c}{2}$, $H^\prime_2(x) < 0$
and when $\dfrac{c}{2} < x \le c - \dfrac{1}{2}$, $H^\prime_2(x) > 0$.

%%%%%%%%%%%%%%%%%%%%%%%%%%%%%%%%%%%%%%%%%%%%%%%%%%%%%%%%%%%%%%%%%%%%%%%%%%%%%%%%%%%%%%%%%%%%%%%%%%%%%%%%%%%%%%%%%%%%%%%%%%%%%%%%%%%%%%%%%%%%%%%%%%%%%%%%
%%%%%%%%%%%%%%%%%%%%%%%%%%%%%%%%%%%%%%%%%%%%%%%%%%%%%%%%%%%%%%%%%%%%%%%%%%%%%%%%%%%%%%%%%%%%%%%%%%%%%%%%%%%%%%%%%%%%%%%%%%%%%%%%%%%%%%%%%%%%%%%%%%%%%%%%

\section{Main Results}

%%%%%%%%%%%%%%%%%%%%%%%%%%%%%%%%%%%%%%%%%%%%%%%%%%%%%%%%%%%%%%%%%%%%%%%%%%%%%%%%%%%%%%%%%%%%%%%%%%%%%%%%%%%%%%%%%%%%%%%%%%%%%%%%%%%%%%%%%%%%%%%%%%%%%%%%
%%%%%%%%%%%%%%%%%%%%%%%%%%%%%%%%%%%%%%%%%%%%%%%%%%%%%%%%%%%%%%%%%%%%%%%%%%%%%%%%%%%%%%%%%%%%%%%%%%%%%%%%%%%%%%%%%%%%%%%%%%%%%%%%%%%%%%%%%%%%%%%%%%%%%%%%

Now, suppose $n > e^{12}$. The product of all primes $p \in (3n, 4n]$, if any, must divide $\binom{4n}{3n}$.
Let $\beta(p)$ be the power of $p$ in the prime factorization of $\binom{4n}{3n}$.
Let
	$$\binom{4n}{3n} = T_1T_2T_3$$
where
	$$T_1 = \prod_{p \le \sqrt{4n}}p^{\beta(p)}, \quad T_2 = \prod_{\sqrt{4n} < p \le 3n }p^{\beta(p)} \quad\text{and}\quad
	  T_3 = \prod_{3n < p \le 4n }p^{\beta(p)}.$$
Bounding each multiplicand in $T_1$ from above by $4n$ (see \cite[p.\ 24]{robbins}) and applying Lemma~\ref{Lemma:oneone},
	$$T_1 < (4n)^{\pi(\sqrt{4n})} \le (4n)^{\frac{\sqrt{4n}}{2}} = (4n)^{\sqrt{n}}.$$
Consider $T_2$. As the prime factorization of $\displaystyle\binom{n}{j}$ in \cite[p.\ 24]{robbins} manifests, for $\sqrt{4n} < p \le 3n$, $\beta(p) \le 1$.

Let $x > 0$ and let $[x]$ be the greatest integer less than or equal to $x$.
Define $\left\{x\right\} = x - [x]$. Let $r$ and $s$ be real numbers satisfying $s > r \ge 1$.
Observe that number of integers in the interval $(s - r, s]$ is $[s] - [s - r]$, which is $[r]$ if $\left\{s\right\} \ge \left\{r\right\}$
and $[r] + 1$ if $\left\{s\right\} < \left\{r\right\}$. Let $N$ be the set of all positive integers.
We define
	$$\begin{Bmatrix}s\\r\end{Bmatrix} = \frac{\displaystyle\prod_{k \in (s - r, s] \cap N} k}{\displaystyle\prod_{k \in (0, r] \cap N} k}
	  = \delta(r, s)\binom{[s]}{[r]},$$
where $\delta(r, s) = 1$ if $\left\{s\right\} \ge \left\{r\right\}$ and $\delta(r, s) = [s - r] + 1$ if $\left\{s\right\} < \left\{r\right\}$.
In both cases, $\delta(r, s) \le s$.

Now let 
$A = \begin{Bmatrix}4n/3\\n\end{Bmatrix}$, 
$B = \begin{Bmatrix}2n\\3n/2\end{Bmatrix}$, 
$C = \begin{Bmatrix}4n/17\\3n/13\end{Bmatrix}$ and
$D = \begin{Bmatrix}2n/7\\4n/15\end{Bmatrix}$.

\medskip

We have the following observations:

\begin{itemize}
	\item $\displaystyle\prod_{\sqrt{4n} < p \le \frac{n}{6}}p \le \prod_{p \le \frac{n}{6}} p \le 4^{\frac{n}{6}}$ (by Lemma~\ref{Lemma:onetwo})
	\item If $\dfrac{n}{6} < p \le \dfrac{2n}{11}$, then
		    $$2p < \frac{n}{2} < 3p < 8p < \frac{3n}{2} < 9p < 11p \le 2n.$$
		    Hence $\displaystyle\prod_{\frac{n}{6} < p \le \frac{2n}{11}}p$ divides $B$.
  \item If $\dfrac{2n}{11} < p \le \dfrac{4n}{21}$, then
		    $$p < \frac{n}{3} < 2p < 5p < n < 6p < 7p \le \frac{4n}{3}.$$
		    Hence $\displaystyle\prod_{\frac{2n}{11} < p \le \frac{4n}{21}}p$ divides $A$. 
	\item If $\dfrac{4n}{21} < p \le \dfrac{n}{5}$, then
		    $$5p \le n < 6p < 15p \le 3n < 16p < 20p \le 4n < 21p.$$
		    Hence $\beta(p) = 0$.
	\item If $\dfrac{n}{5} < p \le \dfrac{2n}{9}$, then
		    $$p < \frac{n}{3} < 2p < 4p < n < 5p < 6p \le \frac{4n}{3}.$$
		    Hence $\displaystyle\prod_{\frac{n}{5} < p \le \frac{2n}{9}}p$ divides $A$.
	\item If $\dfrac{2n}{9} < p \le \dfrac{3n}{13}$, then
		    $$4p < n < 5p < 13p < 3n < 14p < 17p < 4n < 18p.$$
		    Hence $\beta(p) = 0$.
	\item $\displaystyle\prod_{\frac{3n}{13} < p \le \frac{4n}{17}}p$ divides $C$.
	\item If $\dfrac{4n}{17} < p \le \dfrac{n}{4}$, then
		    $$4p \le n < 5p < 12p \le 3n < 13p < 16p \le 4n < 17p.$$
		    Hence $\beta(p) = 0$.
	\item If $\dfrac{n}{4} < p \le \dfrac{4n}{15}$, then
		    $$p < \frac{n}{3} < 2p < 3p < n < 4p < 5p \le \frac{4n}{3}.$$
		    Hence $\displaystyle\prod_{\frac{n}{4} < p \le \frac{4n}{15}}p$ divides $A$.
	\item $\displaystyle\prod_{\frac{4n}{15} < p \le \frac{2n}{7}}p$ divides $D$.
	\item If $\dfrac{2n}{7} < p \le \dfrac{3n}{10}$, then
		    $$3p < n < 4p < 10p \le 3n < 11p < 13p < 4n < 14p.$$
		    Hence $\beta(p) = 0$.
	\item If $\dfrac{3n}{10} < p \le \dfrac{n}{3}$, then
		    $$p < \frac{n}{2} < 2p < 4p < \frac{3n}{2} < 5p < 6p \le 2n.$$
		    Hence $\displaystyle\prod_{\frac{3n}{10} < p \le \frac{n}{3}}p$ divides $B$.
	\item If $\dfrac{n}{3} < p \le \dfrac{4n}{9}$, then
		    $$\frac{n}{3} < p < 2p < n < 3p \le \frac{4n}{3}.$$
		    Hence $\displaystyle\prod_{\frac{n}{3} < p \le \frac{4n}{9}}p$ divides $A$.
  \item If $\dfrac{4n}{9} < p \le \dfrac{n}{2}$, then
		    $$2p \le n < 3p < 6p \le 3n < 7p < 8p \le 4n < 9p.$$
		    Hence $\beta(p) = 0$.
	\item If $\dfrac{n}{2} < p \le \dfrac{2n}{3}$, then
		    $$\frac{n}{2} < p < 2p < \frac{3n}{2} < 3p \le 2n.$$
		    Hence $\displaystyle\prod_{\frac{n}{2} < p \le \frac{2n}{3}}p$ divides $B$.
	\item If $\dfrac{2n}{3} < p \le \dfrac{3n}{4}$, then
		    $$p < n < 2p < 4p \le 3n < 5p < 4n < 6p.$$
		    Hence $\beta(p) = 0$.
  \item If $\dfrac{3n}{4} < p \le \dfrac{4n}{5}$, then
		    $$\frac{n}{2} < p < \frac{3n}{2} < 2p \le 2n.$$
		    Hence $\displaystyle\prod_{\frac{3n}{4} < p \le \frac{4n}{5}}p$ divides $B$.
	\item If $\dfrac{4n}{5} < p \le n$, then
		    $$p \le n < 2p < 3p \le 3n < 4p \le 4n < 5p.$$
		    Hence $\beta(p) = 0$.
	\item $\displaystyle\prod_{n < p \le \frac{4n}{3}}p$ divides $A$.
	\item If $\dfrac{4n}{3} < p \le \frac{3n}{2}$, then
		    $$n < p < 2p \le 3n < 4n < 3p.$$
		    Hence $\beta(p) = 0$.
	\item $\displaystyle\prod_{\frac{3n}{2} < p \le 2n}p$ divides $B$.
	\item If $2n < p \le 3n$, then
		    $$n < p \le 3n < 4n < 2p.$$
		    Hence $\beta(p) = 0$.
\end{itemize}
Therefore, to summarize, we get
	$$T_2 \le 4^{\frac{n}{6}}ABCD.$$
Note that by Lemma~\ref{Lemma:onethree},
 \begin{align*}
   \binom{4n}{3n} &= \frac{(4n)!}{(3n)! n!}\\
                  &> \frac{g(4n)}{f(3n)f(n)}\\
                  &= \frac{2}{\sqrt{6\pi n}}e^{\frac{1}{48n + 1} - \frac{1}{36n} - \frac{1}{12n}}\biggl(\frac{256}{27}\biggr)^n,
 \end{align*}
and similarly,
 \begin{align*}
   A &=   \begin{Bmatrix}4n/3\\n\end{Bmatrix} \le \frac{4n}{3}\binom{[4n/3]}{n} 
      =   \frac{4n}{3}\cdot\frac{\left[\dfrac{4n}{3}\right]!}{n!\biggl(\left[\dfrac{4n}{3}\right] - n\biggr)!}\\
     &<   \frac{4n}{3}\cdot\frac{f\biggl(\left[\dfrac{4n}{3}\right]\biggr)}{g(n)g\biggl(\left[\dfrac{4n}{3}\right] - n\biggr)}\\
     &\le \frac{4n}{3}\cdot\frac{f\biggl(\dfrac{4n}{3}\biggr)}{g(n)g\biggl(\dfrac{n}{3}\biggr)}\quad\text{(by Lemma~\ref{Lemma:onefour})}\\
     &=   \frac{4n}{3}\sqrt{\frac{2}{\pi n}}e^{\frac{1}{16n} - \frac{1}{12n + 1} - \frac{1}{4n + 1}}\biggl(\frac{4^{\frac{4}{3}}}{3}\biggr)^n, 
 \end{align*}
 \begin{align*}
   B &=   \begin{Bmatrix}2n\\3n/2\end{Bmatrix} \le 2n\binom{2n}{[3n/2]} 
      =   2n\cdot\frac{\left[\dfrac{3n}{2}\right] + 1}{2n - \left[\dfrac{3n}{2}\right]}\binom{2n}{[3n/2] + 1}\\
     &<   2n\cdot\frac{\dfrac{3n}{2} + 1}{2n - \dfrac{3n}{2}}\cdot
          \frac{f(2n)}{g\biggl(\left[\dfrac{3n}{2}\right] + 1\biggr)g\biggl(2n - \biggl(\left[\dfrac{3n}{2}\right] + 1\biggr)\biggr)}\\
     &<   (6n + 4)\cdot\frac{f(2n)}{g\biggl(\dfrac{3n}{2}\biggr)g\biggl(2n - \dfrac{3n}{2}\biggr)}\quad\text{(by Lemma~\ref{Lemma:onefive})}\\
     &=   \frac{12n + 8}{\sqrt{3\pi n}}e^{\frac{1}{24n} - \frac{1}{18n + 1} - \frac{1}{6n + 1}}\biggl(\frac{16}{3^{\frac{3}{2}}}\biggr)^n, 
 \end{align*}
 \begin{align*}
   C &=   \begin{Bmatrix}4n/17\\3n/13\end{Bmatrix} \le \frac{4n}{17}\binom{[4n/17]}{[3n/13]} 
      =   \frac{4n}{17}\cdot\frac{\left[\dfrac{3n}{13}\right] + 1}{\left[\dfrac{4n}{17}\right] - \left[\dfrac{3n}{13}\right]}
          \binom{[4n/17]}{[3n/13] + 1}\\
     &\le \frac{4n}{17}\cdot\frac{\dfrac{3n}{13} + 1}{\dfrac{4n}{17} - 1 - \dfrac{3n}{13}}\cdot
          \frac{f\biggl(\left[\dfrac{4n}{17}\right]\biggr)}{g\biggl(\left[\dfrac{3n}{13}\right] + 1\biggr)
          g\biggl(\left[\dfrac{4n}{17}\right] - \biggl(\left[\dfrac{3n}{13}\right] + 1\biggr)\biggr)}
 \end{align*}
 \begin{align*}
    &<   \frac{4n}{17}\cdot\frac{51n + 221}{n - 221}\frac{f\biggl(\dfrac{4n}{17}\biggr)}{g\biggl(\dfrac{3n}{13}\biggr)
          g\biggl(\dfrac{4n}{17} - \dfrac{3n}{13}\biggr)}\quad\text{(by Lemmas~\ref{Lemma:onefour} and \ref{Lemma:onefive})}\\
     &=   \frac{4n}{17}\cdot\frac{51n + 221}{n - 221}\cdot\frac{26}{\sqrt{6\pi n}}
          e^{\frac{17}{48n} - \frac{13}{36n + 13} - \frac{221}{12n + 221}}
          \biggl(221^{\frac{1}{221}}\biggl(\frac{13}{3}\biggr)^{\frac{3}{13}}\biggl(\frac{4}{17}\biggr)^{\frac{4}{17}}\biggr)^n,\\ 
 \end{align*}
and
 \begin{align*}
   D &=   \begin{Bmatrix}2n/7\\4n/15\end{Bmatrix} \le \frac{2n}{7}\binom{[2n/7]}{[4n/15]} 
      =   \frac{2n}{7}\cdot\frac{\left[\dfrac{4n}{15}\right] + 1}{\left[\dfrac{2n}{7}\right] - \left[\dfrac{4n}{15}\right]}
          \binom{[2n/7]}{[4n/15] + 1}\\
     &\le \frac{2n}{7}\cdot\frac{\dfrac{4n}{15} + 1}{\dfrac{2n}{7} - 1 - \dfrac{4n}{15}}\cdot
          \frac{f\biggl(\left[\dfrac{2n}{7}\right]\biggr)}{g\biggl(\left[\dfrac{4n}{15}\right] + 1\biggr)
          g\biggl(\left[\dfrac{2n}{7}\right] - \biggl(\left[\dfrac{4n}{15}\right] + 1\biggr)\biggr)}\\
     &<   \frac{2n}{7}\cdot\frac{28n + 105}{2n - 105}\cdot\frac{f\biggl(\dfrac{2n}{7}\biggr)}{g\biggl(\dfrac{4n}{15}\biggr)
          g\biggl(\dfrac{2n}{7} - \dfrac{4n}{15}\biggr)}\quad\text{(by Lemmas~\ref{Lemma:onefour} and \ref{Lemma:onefive})}\\
     &=   \frac{4n^2 + 15n}{2n - 105}\cdot\frac{15}{\sqrt{2\pi n}}
          e^{\frac{7}{24n} - \frac{5}{16n + 5} - \frac{35}{8n + 35}}
          \biggl(\biggl(\frac{105}{2}\biggr)^{\frac{2}{105}}
          \biggl(\frac{15}{4}\biggr)^{\frac{4}{15}}\biggl(\frac{2}{7}\biggr)^{\frac{2}{7}}\biggr)^n.\\
 \end{align*}
Therefore
	\begin{align*}
		T_3 &= \binom{4n}{3n}\frac{1}{T_1T_2} > \binom{4n}{3n}\frac{1}{(4n)^{\sqrt n}4^{\frac{n}{6}}ABCD}\\
		    &> \frac{\sqrt 3 \pi^{\frac{3}{2}}}{4160}e^EM^n (4n)^{-\sqrt n}\cdot \frac{n^{-\frac{3}{2}}(n - 221)(2n - 105)}{(3n + 2)(3n + 13)(4n + 15)}\\
		    &> \frac{\sqrt 3 \pi^{\frac{3}{2}}}{4160}e^EM^n (4n)^{-\sqrt n}\cdot \frac{n^{-\frac{3}{2}}n^2}{(4n)(4n)(5n)}\\
		    &= \frac{\sqrt 3 \pi^{\frac{3}{2}}}{332800}e^EM^n (4n)^{-\sqrt n} n^{-\frac{5}{2}}\\
	\end{align*}
where
  \begin{align*}
  	E &= \frac{1}{48n + 1} - \frac{1}{36n} - \frac{1}{12n} - \frac{1}{16n} + \frac{1}{12n + 1} + \frac{1}{4n + 1} - \frac{1}{24n} + \frac{1}{18n + 1}
  	     \\
  	  &\hphantom{xx}   + \frac{1}{6n + 1}- \frac{17}{48n} + \frac{13}{36n + 13} + \frac{221}{12n + 221} - \frac{7}{24n} + \frac{5}{16n + 5} + \frac{35}{8n + 35}
  \end{align*}
and
	$$M = \frac{256}{7}\biggl(\frac{1}{4}\biggr)^{\frac{4}{3}}(3)\frac{3^{\frac{3}{2}}}{16}\biggl(\frac{1}{221}\biggr)^{\frac{1}{221}}
				\biggl(\frac{3}{13}\biggr)^{\frac{3}{13}}
	      \biggl(\frac{17}{4}\biggr)^{\frac{4}{17}}\biggl(\frac{2}{105}\biggr)^{\frac{2}{105}}
	      \biggl(\frac{4}{15}\biggr)^{\frac{4}{15}}\biggl(\frac{7}{2}\biggr)^{\frac{2}{7}}4^{-\frac{1}{6}} > 1.$$
Obviously
	$$\lim_{n \rightarrow \infty} e^E = 1.$$
Moreover, we have
	$$\ln \biggl(M^n (4n)^{-\sqrt n} n^{-\frac{5}{2}}\biggr) = n\ln M - \sqrt n \ln (4n) - \frac{5}{2}\ln n.$$
When $n$ tends to infinity, it is easy to check that $\sqrt n \ln (4n) = o(n)$ and $\ln n = o(n)$.
Thus, $\ln \biggl(M^n (4n)^{-\sqrt n} n^{-\frac{5}{2}}\biggr)$ goes to infinity and so does $M^n (4n)^{-\sqrt n} n^{-\frac{5}{2}}$.

It follows that
	$$\lim_{n \rightarrow \infty} T_3 = +\infty,$$
which means that there exists some $n_0$ such that for all $n \ge n_0$, $T_3 > 1$.
In fact, it is routine to check (using WolframAlpha for instance) that when $n > e^{12}$,
$\dfrac{\sqrt 3 \pi^{\frac{3}{2}}}{332800}e^EM^n (4n)^{-\sqrt n} n^{-\frac{5}{2}}$ is always greater than $1$ and so $T_3 > 1$.
Direct verification, on the other hand, ensures that there is always a prime in the interval $[3n, 4n]$ for all positive integers $n < e^{12}$.
Therefore, our desired result ensues:

%%%%%%%%%%%%%%%%%%%%%%%%%%%%%%%%%%%%%%%%%%%%%%%%%%%%%%%%%%%%%%%%%%%%%%%%%%%%%%%%%%%%%%%%%%%%%%%%%%%%%%%%%%%%%%%%%%%%%%%%%%%%%%%%%%%%%%%%%%%%%%%%%%%%%%%%

\begin{Theorem}\label{Theorem:twoone}
For every positive integer $n$, there is a prime in the interval $[3n, 4n]$.
Plainly, it follows that when $n \ge 2$, there is always a prime in the interval $(3n, 4n)$.
\end{Theorem} 

\begin{Corollary}
If $n \ge 3$, then there is a prime in the interval $\biggl(n, \dfrac{4(n + 2)}{3}\biggr)$.
\end{Corollary}

\proof If $n \equiv 0 {\rm~(mod~}3)$, then the result follows directly from Theorem~\ref{Theorem:twoone}.
If $n \equiv 1 {\rm~(mod~}3)$, then by Theorem~\ref{Theorem:twoone} there exists a prime $p \in \biggl(n + 2, \dfrac{4(n + 2)}{3}\biggr)$.
If $n \equiv 2 {\rm~(mod~}3)$, then by Theorem~\ref{Theorem:twoone} there exists a prime $p \in \biggl(n + 1, \dfrac{4(n + 1)}{3}\biggr)$.

\medskip

Next, we establish a lower bound for the number of primes in the interval $[3n, 4n]$.
Bounding each prime in the interval from above by $4n$, we have the following

\begin{Theorem}
For $n \ge 4$, the number of primes in the interval $(3n, 4n)$ is at least
	$$\log_{4n}\biggl(\dfrac{\sqrt 3 \pi^{\frac{3}{2}}}{332800}e^EM^n (4n)^{-\sqrt n} n^{-\frac{5}{2}}\biggr).$$
\end{Theorem}

Note that
	\begin{align*}
		&\log_{4n}\biggl(\dfrac{\sqrt 3 \pi^{\frac{3}{2}}}{332800}e^EM^n (4n)^{-\sqrt n} n^{-\frac{5}{2}}\biggr)\\
		&= \frac{-\dfrac{5}{2}\ln n + n\ln M - \sqrt n \ln (4n) + E + \ln\biggl(\dfrac{\sqrt 3 \pi^{\frac{3}{2}}}{332800}\biggr)}{\ln n + \ln 4}\\
		&= \frac{n\ln M - \sqrt n \ln (4n) + E + \ln\biggl(\dfrac{\sqrt 3 \pi^{\frac{3}{2}}}{332800}\biggr) + \dfrac{5}{2}\ln 4}{\ln n + \ln 4}
		   - \frac{5}{2}\\
		&> \frac{n\biggl(\ln M - \dfrac{\ln (4n)}{\sqrt n}\biggr)}{2\ln n} - \frac{5}{2}.
	\end{align*}
Now check that $\displaystyle \lim_{n \rightarrow \infty} \frac{\ln (4n)}{\sqrt n} = 0$.
Moreover, it is obvious that
	$$\lim_{n \rightarrow \infty} \frac{n}{\ln n} = +\infty.$$
Thus we have the following

\begin{Theorem}
As $n$ tends to infinity, the number of primes in the interval $[3n, 4n]$ goes to infinity.
In other words, for every positive integer $m$, there exists a positive integer $L$ such that for all $n \ge L$,
there are at least $m$ primes in the interval $[3n, 4n]$.
\end{Theorem}

%%%%%%%%%%%%%%%%%%%%%%%%%%%%%%%%%%%%%%%%%%%%%%%%%%%%%%%%%%%%%%%%%%%%%%%%%%%%%%%%%%%%%%%%%%%%%%%%%%%%%%%%%%%%%%%%%%%%%%%%%%%%%%%%%%%%%%%%%%%%%%%%%%%%%%%%
%%%%%%%%%%%%%%%%%%%%%%%%%%%%%%%%%%%%%%%%%%%%%%%%%%%%%%%%%%%%%%%%%%%%%%%%%%%%%%%%%%%%%%%%%%%%%%%%%%%%%%%%%%%%%%%%%%%%%%%%%%%%%%%%%%%%%%%%%%%%%%%%%%%%%%%%

{\bf ACKNOWLEDGEMENTS.} The author is deeply indebted to Dr. Kin Y. Li of the Hong Kong University of Science and Technology and Mr. C.J. Alaban for their generous and invaluable help.

%%%%%%%%%%%%%%%%%%%%%%%%%%%%%%%%%%%%%%%%%%%%%%%%%%%%%%%%%%%%%%%%%%%%%%%%%%%%%%%%%%%%%%%%%%%%%%%%%%%%%%%%%%%%%%%%%%%%%%%%%%%%%%%%%%%%%%%%%%%%%%%%%%%%%%%%
%%%%%%%%%%%%%%%%%%%%%%%%%%%%%%%%%%%%%%%%%%%%%%%%%%%%%%%%%%%%%%%%%%%%%%%%%%%%%%%%%%%%%%%%%%%%%%%%%%%%%%%%%%%%%%%%%%%%%%%%%%%%%%%%%%%%%%%%%%%%%%%%%%%%%%%%

%%%%%%%%%%%%%%%%%%%%%%%%%%%%%%%%%%%%%%%%%%%%%%%%%%%%%%%%%%%%%%%%%%%%%%%%%%%%%%%%%%%%%%%%%%%%%%%%%%%%%%%%%%%%%%%%%%%%%%%%%%%%%%%%%%%%%%%%%%%%%%%%%%%%%%%%
%%%%%%%%%%%%%%%%%%%%%%%%%%%%%%%%%%%%%%%%%%%%%%%%%%%%%%%%%%%%%%%%%%%%%%%%%%%%%%%%%%%%%%%%%%%%%%%%%%%%%%%%%%%%%%%%%%%%%%%%%%%%%%%%%%%%%%%%%%%%%%%%%%%%%%%%

\end{document}